\pdfoutput=1
\documentclass[12pt]{article}
\usepackage{graphics}
\usepackage{color}
\usepackage{epsfig}
\usepackage{latexsym}
\usepackage{amssymb}
\usepackage{multicol}
\usepackage{graphicx}
\usepackage{amsmath}
\usepackage{amsfonts}
\usepackage{enumerate}
\DeclareMathSymbol\square {\mathrel}{AMSa}{"03}
\newtheorem{theorem}{Proposition}

\begin{document}
\title {The Pythagorean Tree: A New Species}
\author{H. Lee Price}
\date{September, 2008}


\maketitle

\begin{abstract} In 1963, the Dutch mathemetician F.J.M. Barning described an infinite, planar, ternary tree whose nodes are just the set of Primitive Pythagorean triples.  Seven years later, A. Hall independently discovered the same tree. Both used the method of uni-modular matrices to transform one triple to another.  A number of rediscoveries have occurred more recently. The purpose of this article is to announce the discovery of an entirely different ternary tree, and to show how it relates to the one found by Barning and Hall.

\end {abstract}

\section{Background}

A {\it Pythagorean triple }is a triple $[a, b, c]$ of positive integers that are the sides 
of a right triangle; in other words, satisfy the Pythagorean equation $a^2+ b^2= c^2$. If $d$  
is the greatest common divisor of $a, b, c$ then $[{\frac{a}{d},\frac{b}{d},\frac{c}{d}]} = [a', b', c']$ is a Pythagorean triple called {\it primitive}, and when $d>1$ the non-primitive triple $[a, b, c]$ equals a 
multiple $d[a', b', c']$ of a primitive triple. Thus the study of 
triples is most often focused on primitive Pythagorean triples (PPTs).

This study has very ancient origins. An extensive list of PPTs is found on a 
cuneiform tablet circa 1800 B.C.E. (see Conway and Guy {\bf[3]} or Robson {\bf[10]}). Van der 
Waerden {\bf[13]} speculates that much earlier knowledge is possible. Though ancient, the field is wide open, judged by the number and depth of recent papers. We mention two central results: {\bf a)} a simple 
two-parameter generation process, and {\bf b)} an infinite ternary tree 
whose set of nodes is just the set of PPTs. The generation by 
parameters is found in Euclid and Diophantus, but the ternary tree was 
overlooked until the second half of the twentieth century. It was first 
discovered in 1963 by F.J.M. Barning {\bf[1]}, and was then independently discovered 
by A. Hall seven years later. A number of rediscoveries have occurred 
more recently.

In this article, we announce the discovery of an entirely different ternary tree. 
Indirectly the new tree grew out of joint work which I and Frank 
Bernhart have been doing for several years (see: {\bf[2]} and {\bf[8]}).

\section{On Parameters, Angles, and Fibonacci Boxes}

\begin{figure}[h]
\centerline{\includegraphics[width=2.77in,height=0.98in]{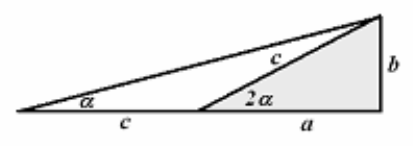}}
\caption{Half-angle tangents}
\label{fig1}
\end{figure}

An easy geometric argument, using figure 1 shows that the 
acute angle marked $2{\alpha} $ in the right triangle with legs 
$({a, b})$ is twice the size of the angle marked $\alpha$ in 
the right triangle with legs $({b, a+c})$, the unshaded triangle being 
isosceles. Hence the tangent of the half-angle $\alpha $ is just ${b}/{(c+a)}={(c-a)}/{b}$. By symmetry, there is another half-angle tangent (or HAT), namely ${a}/{(c + b)}={(c-b)}/{a}$. It 
follows that when the sides $a, b, c$ are rational (or integral) 
then the HATs are rational, and can be reduced to lowest terms.

These two reduced fractions $\frac {q}{p}, \frac {q'}{p'}$ are here displayed as the two columns of a
{$2 \times 2 $} matrix, with the numerators in the first row, and denominators in the second row. The reduced fraction with the largest denominator is put in the second column. These numbers will be revealed below to be the parameters that show up in the two variations of the ancient parametric solution. Some examples are now given.

\[
\mathop {\left[ {{\begin{array}{*{20}c}
 q \hfill & {{q}'} \hfill \\
 p \hfill & {{p}'} \hfill \\
\end{array} }} \right]}\limits_{\left( {\mbox{a,b,c}} \right)} \mbox{  } \rightarrow
\mathop {\left[ {{\begin{array}{*{20}c}
 1 \hfill & 1 \hfill \\
 2 \hfill & 3 \hfill \\
\end{array} }} \right]}\limits_{\left( {\mbox{3,4,5}} \right)} \mbox{ , 
}\mathop {\left[ {{\begin{array}{*{20}c}
 2 \hfill & 1 \hfill \\
 3 \hfill & 5 \hfill \\
\end{array} }} \right]}\limits_{\left( {\mbox{5,12,13}} \right)} \mbox{ , 
}\mathop {\left[ {{\begin{array}{*{20}c}
 1 \hfill & 3 \hfill \\
 4 \hfill & 5 \hfill \\
\end{array} }} \right]}\limits_{\left( {\mbox{15,8,17}} \right)} \mbox{ , 
}\mathop {\left[ {{\begin{array}{*{20}c}
 3 \hfill & 1 \hfill \\
 4 \hfill & 7 \hfill \\
\end{array} }} \right]}\limits_{\left( {\mbox{7,24,25}} \right)} \mbox{ , 
}\mathop {\left[ {{\begin{array}{*{20}c}
 2 \hfill & 3 \hfill \\
 5 \hfill & 7 \hfill \\
\end{array} }} \right]}\limits_{\left( {\mbox{21,20,29}} \right)} 
\]In the third example, $\frac{(17-8)}{15} = \frac{9}{15} = \frac{3}{5}$ and $\frac{(17-15)}{8} = \frac{2}{8} = \frac{1}{4}.$

Several remarkable properties are visible in this display. There is a single 
even number, and it is always in the first column. The column products are 
(almost) the legs of the right triangle. More precisely, the second column 
product equals the odd leg, whereas the first column product is {\it half} the even 
leg. The hypotenuse is the {\it permanent} (the sum of the diagonal products). It is also 
the difference of the two row products:

$${17 = (1\times5) + (3\times4) = (4\times5)-(1\times3)}.$$

Less obviously, the two diagonal products and the two row products are the 
{\it radii} of four circles: the {\bf {\it in-circle}} and the three {\bf {\it 
ex-circles}} of the right triangle. The largest radius (second row product) 
is the sum of the other three. The product of all four numbers is the {\it area} of 
the triangle. In another article {\bf {[2]}} we show that the four radii satisfy the Descartes Circle Equation!

For triple $[5, 12, 13]$ above, the radii are found to be $2, 15$ (row products) and $3, 10$ (diagonal products) 
and the area is  $1 \cdot 2 \cdot 3  \cdot 5 = \frac {1}{2}(5 \cdot 12 ).$  Obviously, $2 \cdot 15 = 3 \cdot 10 = 1 \cdot 2 \cdot 3  \cdot 5.$
Some of these claims are easy corollaries of others. However, the most basic 
fact about the boxes is this: the sum and difference of the first column 
numbers give the second column, that is:  $q' = p - q $ and $p' = p + q.$ These are the \textit{key} equations. Suppose that we convert the $2 \times 2$ into a 4-tuple, $K =[q', q, p, p'],$ by appending the second row to the reversal or flip of the first row. 
The key equations now say that this 4-tuple, $K$ , obeys 
the Fibonacci rule. The sequence $K$ sits in the  
$ 2 \times 2$ box like the letter ``C'', counter-clockwise from the upper right (see fig.2).

\begin{figure}[h]
\centerline{\includegraphics[width=2.16in,height=0.95in]{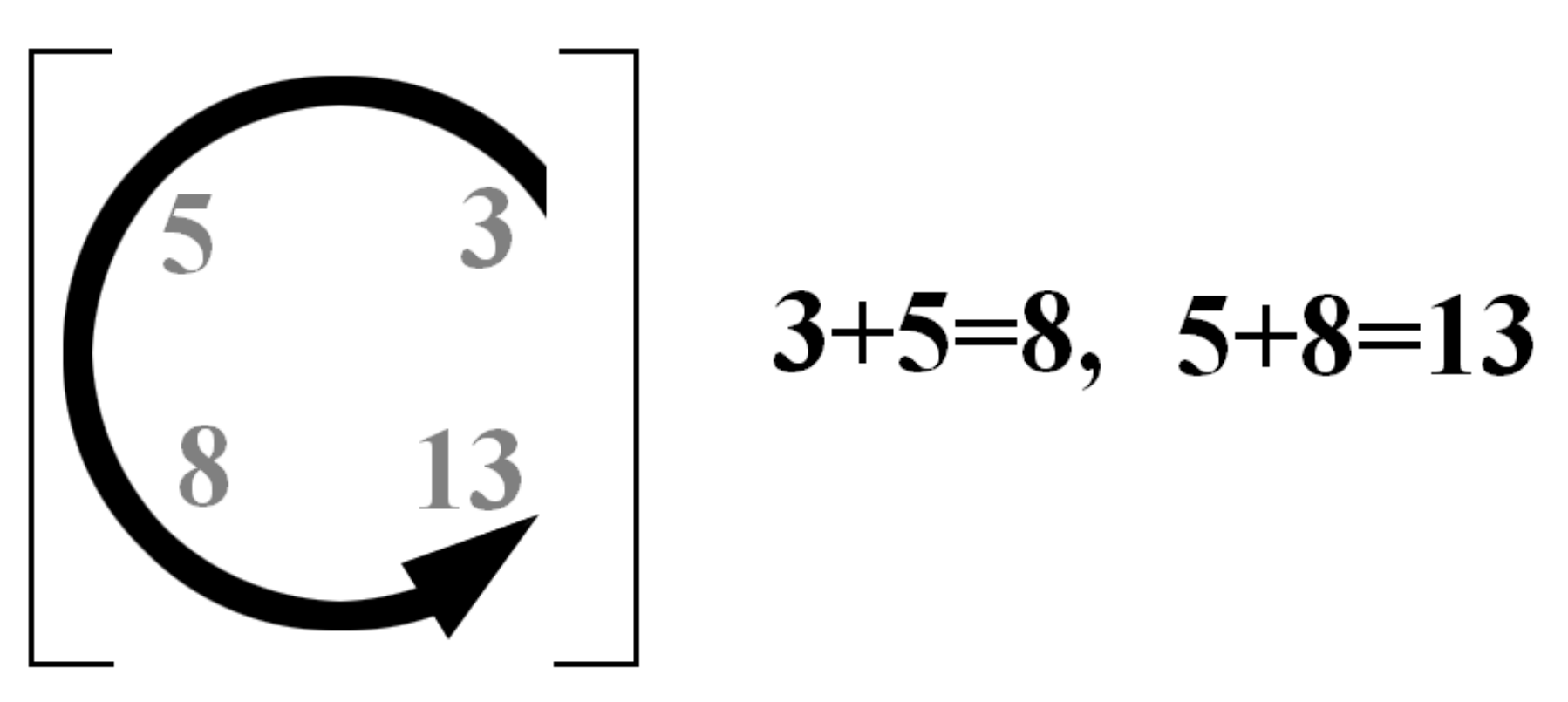}}
\caption{A``Fibonacci box''}
\label{fig2}
\end{figure}

This example, or ``Fibonacci box'', is made up of four consecutive terms of 
the Fibonacci sequence, and it generates the primitive triple $[39, 80, 89].$

The Fibonacci box ({\bf FB}) is uniquely determined by any two numbers in 
it, hence the same is true of the triple. For example pick a column (one of 
the HAT fractions $q/p.)$ The PPT is then $T = [p^{2} - q^{2}, 2pq, p^{2}+ q^{2}]$ if we pick the first column, or $\frac {1}{2}T$ if we pick the second column. These are the two standard 
solutions, known in essence for millennia. We also have a {\it mixed} solution(s) 
employing the column, row, and diagonal products that we gave earlier.
	\[[q'p' ,  2pq,  qp' +q'p ] = [q'p' , 2qp,  pp' - qq'].
\]

We think it is more elegant (and informative!), even though it uses all four parameters, when two 
are sufficient. \newline 
\newline \textbf{Egyptian Fractions and Pythagorean Triples}
\newline

One perhaps surprising application using all four parameters is to the writing of Eqyptian Fractions (EFs). Consider the \textit{perimeter} to \textit{area} ratio given by $\frac{2}{qq'}=\frac{1}{qq'}+\frac{1}{qp'}+\frac{1}{q'p}+\frac{1}{pp'}.$  The sum of the denominators (radii) on the right-hand side \emph{is} the perimeter of the corresponding PPT.   Subtract $\tfrac{1}{qq'}$ from both sides to obtain a 3-term EF variant.  Now, subtract $\tfrac{1}{pp'}$.  The resultant 2-term EF, $\tfrac{1}{q'p}+\tfrac{1}{qp'}$, gives the ratio of \textit{hypotenuse} to \textit{area}.  From the 4-term EF above, we also have $\tfrac{2}{q'p'}=\tfrac{1}{q'p}+\tfrac{1}{pp'}$ which, together with its multiples, produces all of the 2-term solutions found in the Rhind Mathematical Papyrus $2/n$ Table.  Thus each {\bf FB} produces both PPTs and EFs.

In the process of developing our work, Bernhart and I used the key equations 
and the ``HAT'' box for some time before we noticed the Fibonacci property. 
This property is found also in Horodam {\bf [8]}, who cites Dickson {\bf 
[4]} for references to special cases noticed in the 19th century. 
Horodam does not use the same expression(s) for the hypotenuse as we do, and he 
has no reference to the HATs.

\section{Connections}

The reader may wish to fill in the details needed for all the 
claims made in the last section. We give some hints before turning to the 
matter of the trees. The two HATs, $x = \frac {q}{p}$ and $y = \frac {q'}{p'}$, satisfy $xy + x + y = 1.$ This may be proved by substituting the values $x = \frac {(c-a)}{b}$ and $y = \frac {(c-b)}{a}.$

Since the sum of the half angles equals half of a right angle, we could also 
apply the tangent sum formula: $ \tan {45^{\circ}} = 1 = (x+ y)/(1-xy)$ and arrive at the same result. From there it is a small step to the following.

$$\frac {q}{p}= \frac {p' - q'}{p' + q'}, \quad \quad \frac {q'}{p'}= \frac {p - q}{p + q}$$

The right-hand sides are both reduced, up to a factor of two for one case 
only.

A different perspective uses the complex number $z = (p + qi),$ 
where the argument of $z$ is one of the half-angles, and $q/p$ 
is a HAT fraction. Then $z \rightarrow z^{2}= (p^{2} - q^{2}+ 2pqi)$ is a transformation 
that doubles the angle and squares the modulus, obtaining $p^{2} + q^{2}$ for the hypotenuse. Numerous references to 
this approach exist.

A similar technique is to find point $(\frac {a}{c},\frac {b}{c})$ on the unit 
circle $x^{2}+ y^{2}= 1$ in the coordinate plane using a line through $(-1, 0)$ or $(0, -1).$ See Silverman {\bf [12]}, for example. In terms of projective/homogeneous coordinates, ``points'' $(-1, 0, 1)$ and $(0, q, p)$ determine a ``line'' of points $\lambda (-1, 0, 1) + (0, q, p).$ Choose finite $\lambda $ in order to satisfy the Pythagorean equation (this intersects the unit circle at a new point). Rescale to obtain $(p^{2}- q^{2}, 2pq, p^{2}+ q^{2}),$ which must be proportional to $(a, b, c).$

Consider any {\bf FB} that contains four positive integers and satisfies the 
key equations. Define $r_{1},r_{4}$ to be the products in the first and second rows, define $r_{2}, r_{3}$ to be the rising and descending diagonal products, with sum $c,$ and let $b', a$ be the two column products. Put $b = 2b'.$ The product of all four numbers is grouped in three ways to get $ r_{1} r_{4}= r_{2} r_{3}= ab' = ab/2.$

Quite easily $r_{1}+ r_{2}+ r_{3} = r_{4},$ and
$$  a = r_{1}+ r_{2} = r_{4}- r_{3},$$
$$	b = r_{1} + r_{3} = r_{4}- r_{2},$$
$$	c = r_{2}+ r_{3} = r_{4}- r_{1},$$
and $a^{2}+ b^{2}= c^{2}$ is equivalent to $r_{1} r_{4}= r_{2} r_{3},$ which we have just seen. So every {\bf FB} 
gives a Pythagorean triplet, not necessarily primitive.

If some $k >1$ divides all four parameters of the box, divide 
each by $k.$ Clearly this removes a factor of $k^{2}$ 
from the triple.

Now suppose the second column contains an even parameter. Combine key 
equations, obtain $q' + 2q = p'$, and infer both numbers in the second column are even. Divide the 
second column by two, and exchange the two columns. This will result in 
another {\bf FB}, and remove a factor of two from the triple.

When the second column is odd with no common factor $k>1,$ we 
call the box {\it primitive}. In this case the triple is primitive.

When the first row ends in one, the box is clearly primitive. All triples of 
this sort form a family attributed to Pythagoras, also described as the 
triples $[a,b,c]$ with  $b+1 = c$.

If the first row begins with one, and ends with any odd number, again the 
box is clearly primitive. This time the family of primitive triples is 
attributed to Plato, also described as the triples $[a,b,c]$ with $a+2 = c$

We can also get these families in a simple manner, the first by finding 
consecutive square numbers whose difference is an odd square, and the second 
by finding consecutive {\bf {\it odd}} squares whose difference is an even 
square. Both families are central to understanding the Barning-Hall tree \textbf{(BH)}
{\it and} the new tree.

Consider the operation which takes the second row of a primitive box as the 
first row of a new box; i.e. $K \rightarrow K'$ where $K = [\cdots, \cdots,x, y]$ and $K' = [y, x,\cdots,\cdots].$ Then $K'$ is also primitive. For instance
\[
\begin{array}{l}
 \left[ {1,1,2,3} \right] \to \left[ {3,2,5,7} \right] \to  \\ 
 \left[ {7,5,12,17} \right] \to \left[ {17,12,29,41} \right]\begin{array}{*{20}c}
    \to  &  \ldots   \\
\end{array}\left[ {Q_n ,P_n ,P_{n + 1} ,Q_{n + 1} } \right] \\ 
 \end{array}
\]and we encounter the standard Pell sequences $(P_{n}),(Q_{n}).$ Ratios $Q_{n} /P_{n}$ approach $\sqrt{2}$ as a limit. The resulting triples satisfy $\vert a-b \vert = 1$ and were studied by Fermat (see also Eckert {\bf[5]} or 
Hatch {\bf [7]}). This family is simply related to the {\bf BH}-tree. 
It was also studied by Hatch. Recall the two row products 
$(r_{1}, r_{4})$ and the two diagonal products $(r_{2}, r_{3})$ formed from an {\bf FB}. Let $i, 
j, k $ rotate the values $1, 2, 3.$ Then the compact expression $r_{i}+ r_{j} = r_{4} - r_{k}$ finds all three sides of the triangle for us (these three equations were given separately above).

From this it is easy to set up four tangent circles $(K_{i})$ for $i = 1,2,3,4,$ where $K_{i}$ has radius $r_{i}$ and center $C_{i}.$ If we arrange the centers so that $C_{2}, C_{3}$ and $C_{1}, C_{4}$ are the diagonals of a rectangle of dimensions $a \times b,$ the equations show that all six pairs of circles are tangent, with $K_{4}$ enclosing the other three.

We established in {\bf [2]} that four of the six points of tangency lie on a 
line. Reflection in that line exchanges the other two points of tangency. It 
follows that there is a congruent set of circles using the same six points 
of tangency. Comparison of these circles with the four copies of the right 
triangle $(a, b, c)$ contained in the rectangle shows that these 
circles are also the in-circle and ex-circles! So the radii $r_{i}$ are as claimed.

We have now laid the groundwork for viewing the old and the new tree. 

\section{The First Tree}
\begin{figure}[h]
\centerline{\includegraphics[width=4.9583in,height=2.70in]{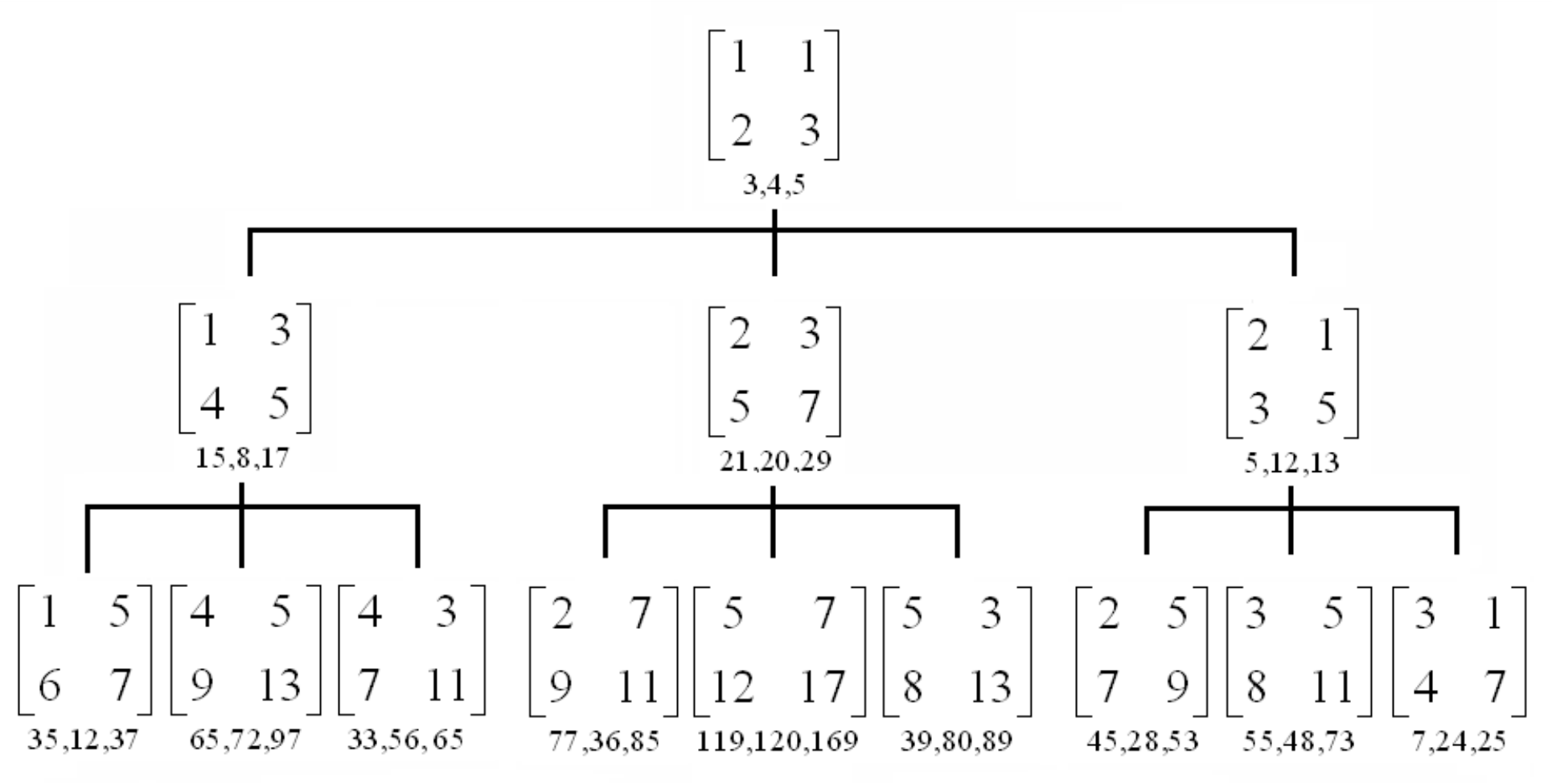}}
\caption{The Barning-Hall Tree (first 3 generations)}
\label{fig3}
\end{figure}

Using our convenient {\bf FB} representation, we can quite easily describe 
both trees. Each one has a first level $(L = 0)$ containing only the 
triple $[3, 4, 5],$ or equivalently, the box containing the 
sequence $K = [1, 1, 2, 3].$

Each node on level $n$  generates three successors (children) on 
the next level. These successors may be called the {\it left child}, {\it middle child}, and {\it right child} (the order is a bit arbitrary). Hence level $L=n$ has $3^{n}$ occupants.

On the Barning-Hall Tree, flipping a column of a box is exchanging the two numbers. If $M$ 
is a primitive Fibonacci box, let $M_{1}, M_{2}, M_{3}$ be obtained from $M$ by flipping the second column, both 
columns, and the first column. The result is 
\underline {not} an {\bf FB}, but there is a unique $M_{i}$ 
which has the {\it same first row}. In other words, discard the second row, and recompute. This 
defines the three children.

It immediately follows from our definition that the first row product of any 
child is the same as one of the other products for the parent. For the 
(left, middle, right) child the first row product is exactly the product in 
the parent of the (main diagonal, second row, back diagonal). Briefly, the 
value $r_{1}$ has been replaced by a different radius $(r_{3}, r_{4}, r_{2}).$ Our name for this is 
\textit{circle promotion}. In effect, one of the three ex-circles is selected, then we say ``you are 
now the in-circle''.

Suppose we have any primitive triple, with its primitive {\bf FB} (box). It 
is straightforward to reverse the process of column flipping, and show that 
this triple/box has a unique place on the tree. Find the positive difference for each row, place the results in the first column (the larger value must go on top), then `complete the box'. Barning and Hall 
both used a different method, the method of matrix transforms. If each  
PPT is written as a column vector, the children are obtained by left 
multiplication by these three matrices. They differ in sign only.

\[
\left[ {\begin{array}{*{20}c}
   { - 1} & 2 & 2  \\
   { - 2} & 1 & 2  \\
   { - 2} & 2 & 3  \\
\end{array}} \right],\quad \left[ {\begin{array}{*{20}c}
   1 & 2 & 2  \\
   2 & 1 & 2  \\
   2 & 2 & 3  \\
\end{array}} \right],\quad \left[ {\begin{array}{*{20}c}
   1 & { - 2} & 2  \\
   2 & { - 1} & 2  \\
   2 & { - 2} & 3  \\
\end{array}} \right]
\]

We omit the details that confirm that our construction and theirs give the 
same tree. Starting at the top, and following \emph{exclusively} the (leftmost, 
middle, rightmost) path downward, one is running through the family of 
triples of (Plato, Fermat, Pythagoras).

For Hall, the tree is the same, except it is rotated ninety degrees, and 
grows from left to right. For Barning the tree is not rotated, but the left 
and middle children are exchanged. Thus the families of Fermat and Plato are 
exchanged.

\section{The Second Tree}
\begin{figure}[h]
\centerline{\includegraphics[width=4.95825in,height=2.81369in]{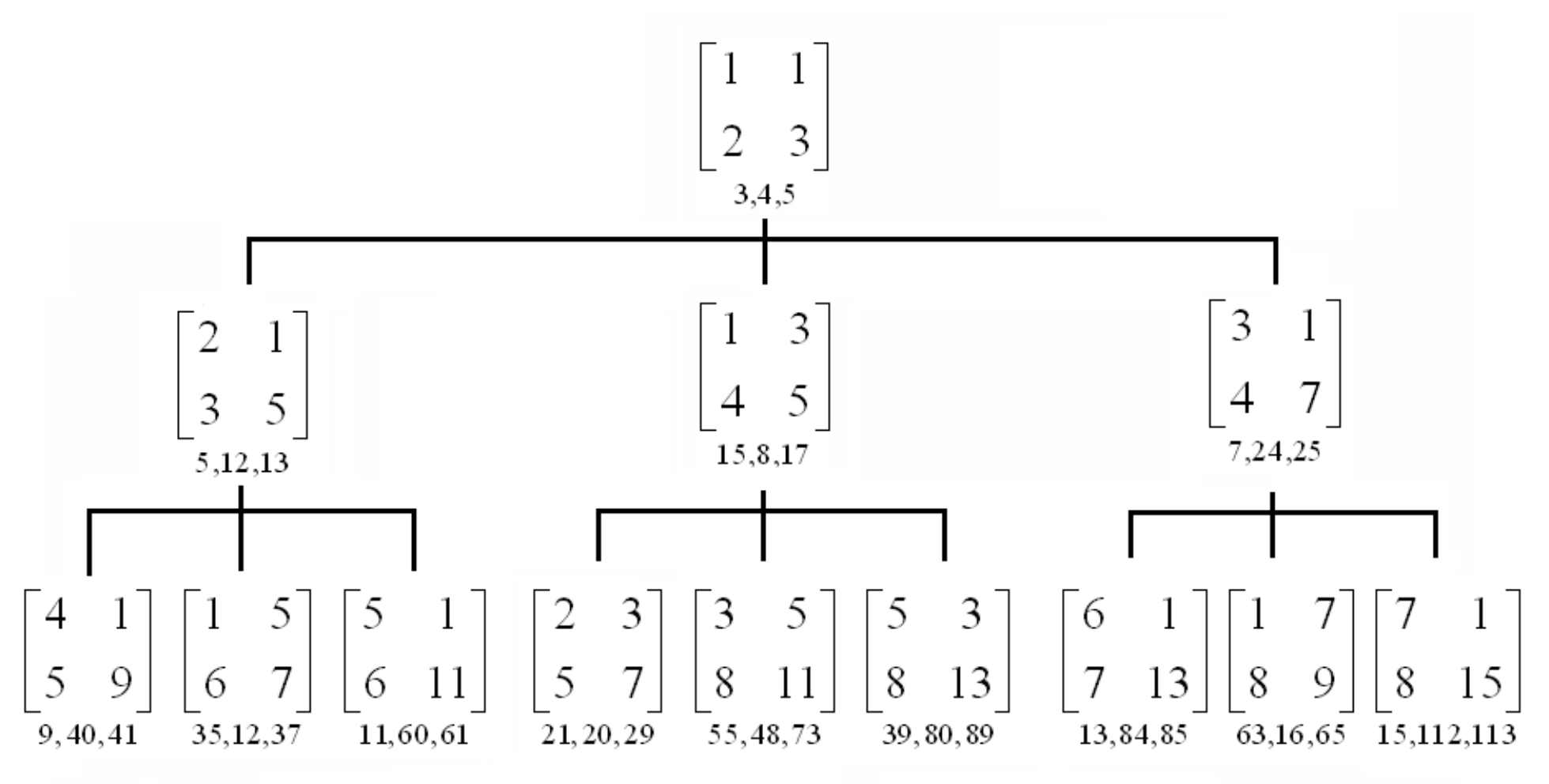}}
\caption{The New Pythagorean Tree (first 3 generations)}
\label{fig6}
\end{figure}
In this section an entirely new tree is revealed. Again we start 
with the box $M,$ which uncoiled, has the sequence $K =[q', q, p, p'.]$ This time we move the last element back one, three, or two steps. That is, $K = [x,\cdots, y,\cdots]$ becomes $K_{1}= [x,\cdots, y,\cdots]$ or $K_{2}= 
[y, x, \cdots,\cdots]$ or $K_{3} = [x, y,\cdots,\cdots]$
Note the slight change of order.

Using these templates, we now construct valid boxes $K_{i}$ such 
that the second column elements $x, y$ of $M$ are 
repositioned as indicated. The first column elements $a, b$ do get 
altered. A visual display looks like this.

\[
M = \left[ {\begin{array}{*{20}c}
   a & x  \\
   b & y  \\
\end{array}} \right]\begin{array}{*{20}c}
   {} & {}  \\
\end{array}K_1  = \left[ {\begin{array}{*{20}c}
   {2a} & x  \\
   y &  \bullet   \\
\end{array}} \right]\begin{array}{*{20}c}
   {} & {}  \\
\end{array}K_2  = \left[ {\begin{array}{*{20}c}
   {x} & y  \\
   {2b} &  \bullet   \\
\end{array}} \right]\begin{array}{*{20}c}
   {} & {}  \\
\end{array}K_3  = \left[ {\begin{array}{*{20}c}
   y & x  \\
   {2b} &  \bullet   \\
\end{array}} \right]
\]

Notice here that one element of the first column is bumped, and the other is 
{\it doubled}. 

The reason that this works is that there is a simple inverse to the 
operation. It goes as follows. Move the odd element of the first column into 
the bottom right corner, flip the second column if necessary (the smaller 
number must be on top), divide the even element of the first column by two, 
and recompute the missing element. Here is an example that reduces 
completely in six steps. The details of the first step are spelled out, then 
all six steps are summarized (see figure 4).

\begin{figure}[h]
\centerline{\includegraphics[width=2.976in,height=1.96in]{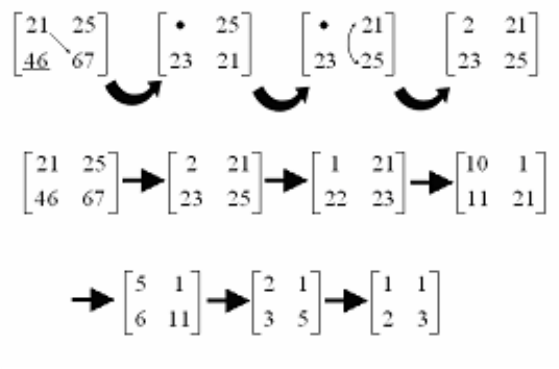}}
\caption{Reversal process in six steps}
\label{fig4}
\end{figure}

The process must end because the final element is always decreasing. In our 
example the final box contains $K_{0}= [1,1,2,3]$ 
and represents the simplest triple $T_{0} = [3,4,5].$ But the reverse process in general stops only because the first 
row is constant, which is $K = [k, k,\cdots,\cdots] = [k, k, 2k, 3k] = 
k[1, 1, 2, 3] = kK_{0}.$

The first column of the box contains the HAT fraction $q/p.$ 
We recall that this is the traditional generator of the triple: $a= p^{2} - q^{2}, b= 2pq, c= p^{2}+ q^{2}.$ There is a simple description of the change in this fraction when we go one step up the tree.

Divide the number that is even by two, then replace the other number by the 
positive difference. For example 

$$
\frac {13}{24} \rightarrow \frac {(13-12)}{12}= \frac {1}{12} \rightarrow \frac {(6-1)}{6}
= \frac {5}{6} \rightarrow \frac {(5-3)}{3} = \frac {2}{3}\rightarrow \frac {1}{(3-1)}= \frac {1}{2}.
$$

Again, if a triple is written as a column vector $v,$ there are 
matrices $M_{1}, M_{2}, M_{3}$ such 
that the three successors of $v$ are $ M_{i}v, i = 1,2,3.$ The determinants are $8, -8, -8.$

\begin{figure}[h]
\centerline{\includegraphics[width=4.472in,height=0.72in]{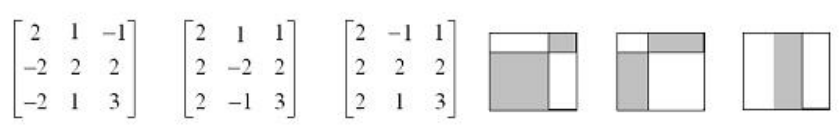}}
\caption{Matrices $M_1 ,M_2,$ and $M_3$.}
\label{fig5}
\end{figure}

The three matrices of Figure 5 are alike, except for signs. The small boxes with 
shading at the right show where the sign differences are between each pair 
of matrices (first two, first and last, last two). The New Tree is shown in Figure 6 below.
The correctness of the matrices may be assured by checking that they 
transform as predicted the three triples $[3,4,5],$ $[5,12,13],$ $[15,8,17],$ since as vectors these three are 
independent. It is also easy to verify that when you transform any triple 
$[a, b, c],$ you are taking one PPT to another PPT. 
E.g. the following simplifies to $a^{2} + b^{2} = c^{2}.$

$$(2a - b + c)^{2} + (2a + 2b + 2c)^{2} =(2a + b + 3c)^{2}$$Now recall the two families of Plato and Pythagoras, where:

\[
\begin{array}{l}
 A(1) = B(1) = \left[ {3,4,5} \right] \\ 
 A(n) = \left[ {4n^2  - 1,} \right],4n,\left[ {4n^2  + 1} \right] \\ 
 B(n) = \left[ {2n + 1,2n(n + 1),2n(n + 1) + 1} \right] = \left[ {a,b,b + 1} \right] \\ 
 \end{array}
\]In the new tree, the three successors of $B(n)$ are $B(2n), A(n+1),$ and $B(2n+1).$ This is easily verified 
by applying the backward step to the HATs $\frac {2n}{2n+1}, \frac{2n+1}{2n+2}, \frac{1}{2n+2}.$ We summarize our results and add a few easy corollaries in the following proposition.

\begin {theorem} 
\item{(a)} The set of PPTs forms an infinite tree with $[3,4,5]$ at 
the root, and the successors of (column vector) $v$ are 
$M_{1} v, M_{2} v, M_{3} v ,$ ($M_{i}$ given above).

\item{(b)} Let the primary HAT (first column of the {\bf FB} of any 
triple $v$ be written as reduced positive fraction $\frac{q}{p}.$ For the successors this becomes 
$\frac {2q}{p+q},\frac {p-q}{2p},\frac {p+q}{2p}.$ Thus the predecessor is 
found by dividing the even number by two and changing the other number to 
the positive difference.

\item{(c)} From the infinite ternary tree, extract an infinite 
binary tree by keeping only the first and third successors. The binary tree 
comprises the family $B(n)$ where $n= 1,2,3,{\ldots}$ reads from left to right, 
level by level.

\item{(d)} The middle successor of $B(n)$ is $A(n).$ So to reach any triple $v$ not in either family, one must start with 
$B(1) = [3,4,5]$ and proceed to $B(n)$ by $k$ steps, where the base two numeral for $n$ has 
$k+1$ digits, thence to $A(n+1),$ and further steps.

\end{theorem}

Part $(d)$ of the {\it Proposition} reads like the philosophical maxim: ``Plato follows Pythagoras, 
and all else follows Plato.'' \newpage

\begin{flushleft}{\textbf{Fermat's Famous Triple}}
\end{flushleft}Any location on the tree can be reached by a path from the top, consisting of $n$ steps to reach the $n^{th}$  level. One step down-left (A), straight-down (B), or down-right (C).   From Figure 6 we see that path code CCC takes us to triple $[15,112,113].$  Path code
\begin{center}
\small{AABAABACABACABAAABABABAAAAAABAAA}
\end{center}
takes us to a famous PPT studied by Fermat. According to Sierpinski {\bf[11]},  Fermat wrote a letter to Mersenne in the year 1643, in which he asserted that primitive triple $$[ 4,565,486,027,761; \quad    1,061,652,293,520;  \quad   4,687,298,610,289 ]$$is the \emph{smallest} in which both the hypotenuse and the leg-sum are perfect squares.  In another article \textbf{[9]} we show that the route to Fermat's triple on the \textbf{BH}-tree is even longer (43 levels down compared to 32 on the New Tree).  The route required is
\begin{center}
\small{BCCCBAAAAAAAAACAABCCCCCCCCCCCCCCCCBCCBAAA}
\end{center}Thus Fermat's triple is one of $3^{43}$ (on level 43, \textbf{BH}-tree)  compared to one of $3^{32}$ ( level 32 on the New Tree).

\section{Concluding Remarks}
An inexact general contrast between the two trees is this: the Barning-Hall tree 
is closely related to the geometry, and the circles, whereas the new tree is 
seemingly less geometric and more number theoretic.
That there is more than one tree constructed by intrinsic means is somewhat 
surprising. Arbitrary arrangements can be concocted, but what counts here is 
the use of a set of three matrices, composed of constants, to generate the 
three successors of a fixed triple. There does not seem to be such a tree in 
which the Pythagorean family $B(n)$ are supplanted by the 
Platonic triples $A(n).$ Frank Bernhart has found a binary tree for 
$A(n)$ corresponding completely to my binary subtree, complete 
with two matrices to generate the two successors. Application of the 
matrices to arbitrary PPTs does not give PPTs.

\underline{\bf AMS Classifications}
\begin{itemize}
\item 11-99  Number theory
\item 11D68 Diophantine equations
\item 11B39 Fibonacci and Lucas numbers
\item 05-B07 Triple systems
\item 05C50 Graph theory
\end{itemize}

\underline{\bf Author Contact Information}

\small{H. Lee Price, 83 Wheatstone Circle, Fairport, NY, USA 14450-1138.}

email: \textit{tanutuva@rochester.rr.com} \newline

[{\LaTeX} document]

\end{document}